# Convergence of impact measures and impact bundles

Leo Egghe, Hasselt University, Belgium

[leo.egghe@uhasselt.be](leo.egghe@uhasselt.be)

ORCID: 0000-0001-8419-2932

**Abstract**

**Purpose** – A new point of view in the study of impact is introduced.

**Design/methodology/approach** – Using fundamental theorems in real analysis we study the convergence of well-known impact measures.

**Findings** – We show that pointwise convergence is maintained by all well-known impact bundles (such as the h-, g-, and R-bundle) and that the µ-bundle even maintains uniform convergence. Based on these results, a classification of impact bundles is given.

**Research limitations** – As for all impact studies, it is just impossible to study all measures in depth.

**Practical implications** – It is proposed to include convergence properties in the study of impact measures.

**Originality/value** - This article is the first to present a bundle classification based on convergence properties of impact bundles.

## Introduction

We denote by $U$ the set of all continuous, decreasing functions defined on the interval $[0,T]$, $T > 0$, with values in $R^+$, where $R^+$ denotes the non-negative real numbers. Note that a function $Z \in U$, does not have to be strictly decreasing and hence $U$ contains all constant functions, including the zero function $O$. As in previous investigations, see e.g. (Egghe, 2021; Egghe & Rousseau, 2022) the functions in $U$ are continuous models for general rank-frequency functions such as authors and their articles (ranked in decreasing order of their numbers of publications); articles co-authored by one scientist and the received number of citations, and so on (Rousseau et al., 2018, p. 301).

Let m be any function from $U$ to $R^+$, thought of as a being a bibliometric measure, such as 'total number of citations', but for the moment without any special property. The first problem we want to study is the following:

If the sequence $(Z_n)_n$ tends to Z does this imply that $m(Z_n)$ tends to $m(Z)$?

Of course, this question must be made more specific, in particular, what do we mean by the expression $(Z_n)_n$ tends to Z, with all $Z_n$ in **U**? Does Z have to belong to **U** too? This will be explained in the next section.

The intuitive meaning of the convergence question is obvious. If rank-frequency functions $Z_n$, e.g., received citations of publications, are given, and if these are evaluated by a measure m, then one only wants to use measures m such that if the sequence $(Z_n)_n$ tends to Z then $m(Z_n)$ tends to $m(Z)$. Stated more loosely: if two cases are very similar, then functions measuring their impact must result in values that do not differ much.

## Convergence of impact measures

Definition. Pointwise convergence

We say that $(Z_n)_n \rightarrow Z$, pointwise, with all $Z_n$ in **U** iff

$$\forall x \in [0,T] : \lim_{n \to \infty} Z_n(x) = Z(x) \text{ in } \mathbf{R}^+.$$

Here we can make a distinction between the case that $Z \in \mathbf{U}$ (pointwise convergence in **U**) and the case that $Z \notin \mathbf{U}$. In the latter case we will say that there is pointwise convergence on **U**.

Definition: uniform convergence

We say that $(Z_n)_n \rightarrow Z$, uniformly in **U** iff

$\forall \varepsilon > 0$, $\exists n_0$ such that $\forall x \in [0,T]$: $n \geq n_0 \Rightarrow |Z_n(x) - Z(x)| < \varepsilon$,

with $Z \in \boldsymbol{U}$

The point is that $n_0$ does not depend on x. We further note that the uniform limit of continuous functions is continuous (Apostol, 1967, 11.3). It is obvious that when $(Z_n)_n \to Z$, uniformly in $\boldsymbol{U}$ then also $(Z_n)_n \to Z$, pointwise in $\boldsymbol{U}$.

## Mathematical preliminaries

This section is purely mathematical. It can be skipped by those readers who are only interested in the informetric applications.

We first recall the following theorem from advanced calculus. Here, and everywhere else integrals are Riemann integrals, not the more general Lebesgue integrals.

Theorem 1 (Apostol, 1967, 11.4)

If $(f_n)_n \to f$, uniformly in $\boldsymbol{U}$ then, for all x in [0,T]: $\int_0^x f_n(s)ds \to \int_0^x f(s)ds$.

It can be shown (Apostol, 1967) that this theorem does not hold for pointwise convergence.

In the next step, we will apply the theorem of Arzelà (Apostol, 1957, p. 405). This theorem states that if a sequence of real, integrable functions $(f_n)_n$, defined on the finite interval [a,b], is uniformly bounded, converging pointwise to a function f, which is integrable on the interval [a,b], then $\int_a^b f_n(s)ds \to \int_a^b f(s)ds$.

We recall that these functions are uniformly bounded if there exists $M \geq 0$ such that for all n and all $x \in [a,b]$: $|f_n(x)| \leq M$.

Arzelà's theorem leads to the following result.

Theorem 2.

If the decreasing functions $(Z_n)_n \to Z$, pointwise, with all $Z_n$ in **U**, then $\int_0^T Z_n(s)ds \to \int_0^T Z(s)ds$.

Proof.

Applying Arzelà's theorem to our situation shows that if the decreasing functions $(Z_n)_n$ converge pointwise to Z, then Z may not be continuous, but it is certainly integrable as a function of bounded variation. Moreover, the sequence $(Z_n)_n$ is uniformly bounded, because $(Z_n(0))_n$ is bounded, hence we have that $\int_0^T Z_n(s)ds \to \int_0^T Z(s)ds$.

Remark 1. We did not have to assume in this application that the limiting function Z is continuous.

Remark 2. Arzelà's theorem does not imply that the pointwise convergence of the sequence of real, integrable functions $(f_n)_n$, is automatically uniform. Indeed, the following example gives a sequence of, uniformly bounded, continuous functions pointwise converging to a discontinuous function f, proving that the convergence cannot be uniform (recall that in the case of uniform convergence the limiting function of continuous functions is always continuous). We define $f_n(x) = 1 - x^n$ on the interval [0,1]. Then $f(x) = 1$ for x in [0,1[, while $f(1) = 0$.

Next, we mention for further use, the so-called second theorem of Dini (but probably due to Pólya).

Theorem: Second Dini theorem

Let $(f_n)_n$ be a sequence of increasing or decreasing real functions, defined on [a,b], such that $(f_n)_n$ tends to f, pointwise and f is continuous, then $(f_n)_n$ tends to f uniformly.

Note that the functions $f_n$ in Dini's second theorem do not have to be continuous.

For the mathematically inclined readers, we note that we assume that we work in first-countable topological function spaces so that using sequences suffices (Kelley, 1975).

## Informetric applications

*Total and average number of items on an interval starting in 0*

Consider $\theta \in$ [0,T] fixed, and define $I_\theta(Z) = \int_0^\theta Z(s)ds$ and $\mu_\theta(Z) = \frac{1}{\theta}\int_0^\theta Z(s)ds$, where $\mu_0(Z)$ is defined as a limit and is equal to Z(0).

If Z represents citations of publications, then the operator $I_\theta$ is the continuous, hence model-theoretic, expression for the total number of citations received by the first $\theta$ publications, where publications are ranked in decreasing order according to received citations. Similarly, $\mu_\theta$ is the continuous expression for the average number of citations received by the first $\theta$ publications.

Then, by the previous results, for fixed $\theta \in\ ]0,T]$, $(Z_n)_n \to Z$, pointwise in $\boldsymbol{U}$ implies that $I_\theta(Z_n) \to I_\theta(Z)$, and similarly $(Z_n)_n \to Z$, pointwise in $\boldsymbol{U}$ implies that $\mu_\theta(Z_n) \to \mu_\theta(Z)$. As uniform convergence implies pointwise convergence, the same results hold for uniform convergence.

*The g-index*

We next show that a similar property holds for the generalized g-index (van Eck & Waltman, 2008). We recall from (Egghe & Rousseau, 2019) that if $Y(x) = \int_0^x Z(s)ds$, with $x \in [0,T]$, $Z \in \boldsymbol{U}$, and $Y(T) \leq \theta\ T^2$, then there exists a unique point $g_\theta$ in $[0,T]$ such that $Y(g_\theta) = \theta(g_\theta)^2$. These $g_\theta$ values are called generalized g-indices. If $\theta = 1$ we obtain the g-index as introduced by Egghe (2006a,b).

Theorem 3.

Let $(Z_n)_n \to Z$, pointwise on $\boldsymbol{U}$, and let $\theta \in \boldsymbol{R}^+$ such that $g_\theta(Z)$ and $g_\theta(Z_n)$ exist for every n, then $g_\theta(Z_n) \to g_\theta(Z)$, with $\theta$ fixed.

Proof. Take $\theta > 0$, fixed, and assume that $(Z_n)_n \to Z$, pointwise. By the definition of the generalized g-index we have:

$x_n = g_\theta(Z_n) \Leftrightarrow \int_0^{x_n} Z_n(s)ds = \theta\ x_n^2$   and   similarly:   $x = g_\theta(Z) \Leftrightarrow \int_0^x Z(s)ds = \theta\ x^2$

Hence: $|x_n^2 - x^2| = \frac{1}{\theta}\left|\int_0^{x_n} Z_n(s)ds - \int_0^x Z(s)ds\ \right|$

$$\frac{1}{\theta}\left|\int_0^x Z_n(s)ds + \int_x^{x_n} Z_n(s)ds - \int_0^x Z(s)ds\right|$$

$\leq \frac{1}{\theta}\left(\left|\int_0^x Z_n(s)ds - \int_0^x Z(s)ds\right| + \left|\int_x^{x_n} Z_n(s)ds\right|\right).$

Now: $\left|\int_x^{x_n} Z_n(s)ds\right| = |x_n - x| \times$ (the average value of $Z_n$ on the interval $[x, x_n]$ (or $[x_n, x]$, whichever applies)). As the functions $Z_n$ are decreasing this average is smaller than or equal to the average taken over $[0, x]$. This is, by definition, $|x_n - x| * \theta x_n$.

Hence, $\quad |x_n^2 - x^2| = |x_n - x| * (x_n + x) \leq \frac{1}{\theta}\left|\int_0^x Z_n(s)ds - \int_0^x Z(s)ds\right| \quad +$ $|x_n - x| * x_n$. Then:

$$|x_n - x|x \leq \frac{1}{\theta}\left|\int_0^x Z_n(s)ds - \int_0^x Z(s)ds\right|$$

or $\quad |g_\theta(Z_n) - g_\theta(Z)| = |x_n - x| \leq \frac{1}{\theta x}\left|\int_0^x Z_n(s)ds - \int_0^x Z(s)ds\right| \quad$ . By Arzelà's theorem applied on the interval $[0,x]=[0,g_\theta(Z)]$, this proves Theorem 3.

Corollary

Let $(Z_n)_n \to Z$, uniform in **U** and let $\theta \in \mathbf{R}^+$ such that $g_\theta(Z)$ and $g_\theta(Z_n)$ exist for every n, then $g_\theta(Z_n) \to g_\theta(Z)$, with $\theta$ fixed.

*PED measures and the h-index*

We next study the case of the so-called (PED)-measures, including the $h_\theta$-measures, with $\theta$ fixed. We recall that a (PED)-measure m is a measure such that there exists a continuous strictly increasing function $f_m$ defined on $[0,T]$ such that $x = m(Z)$ iff $Z(x) = f_m(x)$ (Egghe, 2021). Note, that in the case of (PED)-measures, $x = m(Z) \in [0,T]$. The $h_\theta$-measures are PED measures such that $f_m(x) = \theta x$ ($\theta$ fixed).

Theorem 4.

If m is a (PED)-measure, with associated function $f_m$, then $(Z_n)_n$ → Z, pointwise, with all $Z_n$ in **U** implies that $(m(Z_n))_n → m(Z)$.

Proof. We put $y_n = m(Z_n)$ and $y = m(Z)$. Then, by definition

$$y_n = m(Z_n) \Leftrightarrow Z_n(y_n) = f_m(y_n)$$

$$y = m(Z) \Leftrightarrow Z(y) = f_m(y)$$

We observe that these equivalences are consequences of the facts that the sequence $(Z_n)_n$ is decreasing and that $f_m$ is strictly increasing. Now we have, because f is strictly increasing,

$$|f_m(y) - f_m(y_n)| = \begin{cases} f_m(y) - f_m(y_n) & if \ y \geq y_n \\ f_m(y_n) - f_m(y) & if \ y \leq y_n \end{cases}$$

In the first case we have $Z_n(y) \leq Z_n(y_n)$, because each $Z_n$ is decreasing, leading to $0 \leq Z(y) - Z_n(y_n) \leq Z(y) - Z_n(y)$.

In the second case, again using the fact that all $Z_n$ are decreasing, we have $0 \leq Z_n(y_n) - Z(y) \leq Z_n(y) - Z(y)$.

Consequently: $|f_m(y) - f_m(y_n)| = |Z(y) - Z_n(y_n)| \leq |Z(y) - Z_n(y)|$ . Because, now, $(Z_n)_n$ → Z, pointwise, we have $\lim_{n\to\infty} f_m(y_n) = f_m(y)$. We know from real analysis (De Lillo, 1982, p. 119) that if a function f is strictly monotone on a closed interval, then its inverse function $f^{-1}$ exists and is continuous on the interval [f(0), f(T)] (in the case that f is increasing). Applying this result on $f_m$ and the interval [0,T] leads to

$\lim_{n\to\infty} m(Z_n) = \lim_{n\to\infty} y_n = \lim_{n\to\infty} (f_m^{-1}(f_m(y_n))) = f_m^{-1}(f_m(y))$ = y = m(Z) □

Corollary 1

If m is a (PED) measure, with associated function $f_m$, then $(Z_n)_n \to Z$, uniformly, with all $Z_n$ in **U** implies that $(m(Z_n))_n \to m(Z)$.

Corollary 2: the generalized h-index

If $(Z_n)_n \to Z$, pointwise, with all $Z_n$ in **U** then $(h_\theta(Z_n))_n \to h_\theta(Z)$, with θ fixed.

Corollary 3

Besides for the generalized h-indices, Theorem 4 also holds for the generalized Kosmulski-indices $h_\theta^{(p)}$, with θ fixed (Egghe, 2021).

Corollary 4

Let $\rho_X$: φ → $\rho_X(\varphi)$ denote the polar function of X. If $(Z_n)_n \to Z$, pointwise in **U**, then also $\rho_{Z_n} \to \rho_Z$ pointwise.

Note that the polar functions $\rho_Z$ and all $\rho_{Z_n}$ exist, because Z and all $Z_n$ are decreasing continuous functions.

Proof. From (Egghe & Rousseau, 2020) we know that

$$\rho_Z(\varphi) = h_\theta(Z)\sqrt{1 + \theta^2} \text{ , with } \theta = tg(\varphi)$$

and similar expressions for all $Z_n$. From Theorem 4 and the continuity of the tangent function tg, it follows that $(Z_n)_n \to Z$, pointwise in **U** implies that, for all φ ∈ $\left[0, \frac{\pi}{2}\right]$:

$$h_{tg(\varphi)}(Z_n) \to h_{tg(\varphi)}(Z)$$

Recall that $m = h_{tg(\varphi)}$ is a (PED)-measure, with φ fixed, and $f_m(x)$ = x.tg(φ). As θ = tg(φ) is fixed, this leads to

$$\rho_{Z_n}(\varphi) = h_{tg(\varphi)}(Z_n)\sqrt{1 + (tg(\varphi))^2} \to h_{tg(\varphi)}(Z)\sqrt{1 + \big(tg(\varphi)\big)^2} = \rho_Z(\varphi).$$

Hence, we have $\rho_{Z_n} \to \rho_Z$ pointwise.□

Remark 1. We do not know a place in the mathematical literature where this property of polar functions is proved explicitly.

Remark 2. Corollary 4 does not hold for all monotone functions used in informetric studies. We present a simple example. Let $(Z_n)_n$ be a sequence of convex Lorenz curves (hence defined on [0,1]) converging pointwise to the diagonal Z(x) = x. Then all functions $\rho_{Z_n}$ exist, but $\rho_Z$ does not.

*Percentiles*

The next theorem, dealing with percentiles as a measure, is essentially trivial.

Theorem 5

If all $Z_n \in$ **U** and θ ∈ [0,T] is fixed, then the θth percentile of $Z_n$, defined as $P_\theta(Z_n) = Z_n(\theta)$ satisfies $(Z_n)_n \to Z$, pointwise on **U** implies $P_\theta(Z_n) \to P_\theta(Z)$.

Proof. This follows immediately from the fact that for any X in **U**, P(X) = X.

Corollary

If all $Z_n \in$ **U** and θ ∈ [0,T] is fixed, then $(Z_n)_n \to Z$, uniform in **U** implies $P_\theta(Z_n) \to P_\theta(Z)$.

*The R-index*

Finally, we also consider the R-index (Jin et al., 2007). Recall that for X ∈ **U**, and if $h_\theta$ exists, $R_\theta^2(X) = \int_0^{h_\theta(X)} X(s)ds$

Theorem 6

If all $Z_n \in$ **U** and θ ∈ [0,T], fixed, $(Z_n)_n \to Z$, pointwise on **U** implies $R_\theta(Z_n) \to R_\theta(Z)$.

Proof. $\left| R_\theta^2(Z_n) - R_\theta^2(Z) \right| = \left| \int_0^{h_\theta(Z_n)} Z_n(s)ds - \int_0^{h_\theta(Z)} Z(s)ds \right|$

$$= \left| \int_0^{h_\theta(Z)} Z_n(s)ds + \int_{h_\theta(Z)}^{h_\theta(Z_n)} Z_n(s)ds - \int_0^{h_\theta(Z)} Z(s)ds \right|$$

$$\leq \left| \int_0^{h_\theta(Z)} Z_n(s)ds - \int_0^{h_\theta(Z)} Z(s)ds \right| + \left| \int_{h_\theta(Z)}^{h_\theta(Z_n)} Z_n(s)ds \right|$$

Now, when n increases, the first term in this sum converges to zero by Theorem 2 (based on Arzelà's theorem) (and because θ is fixed), while the second term is smaller than $\max_n(Z_n(0), Z(0)) * |h_\theta(Z_n) - h_\theta(Z)|$. The first factor of this second term exists because a pointwise convergent decreasing sequence of functions on **U** is uniformly bounded, while the second factor tends to zero, because of Corollary 2 of Theorem 4.

## Consequences for bundles

In the previous section, we studied the convergence of measures such as $I_\theta$, $\mu_\theta$, $g_\theta$, $P_\theta$, and $R_\theta$ and more generally (PED)-measures, including $h_\theta$, with $\theta$ fixed and admissible (in the sense that the corresponding measure is well-defined). Taking now $\theta$ variable leads to pointwise convergence of bundles.

We recall the following definition, adapted to our needs, from (Egghe & Rousseau, 2022).

Definition: a bundle

A bundle m is a set of functions, referred to as measures, $m_\theta$, with $\theta$ belonging to a subset of $[0, +\infty]$, detailed further on. These measures are defined on a subset $\boldsymbol{Z} \subset \boldsymbol{U}$. For fixed $Z \in \boldsymbol{Z}$ we have a function $\theta \rightarrow m_\theta(Z)$, where now $\theta$ ranges in a subset of $[0, +\infty]$, depending on Z.

When studying two functions Z and Y at the same time, we will always assume that $\theta$ belongs to the set where $m_\theta(Z)$, as well as $m_\theta(Y)$, are defined. We simply write "all admissible $\theta$". We refer the reader to the previous sections for examples of measures $m_\theta$.

Without new proofs necessary, we have the following theorem.

Theorem 7

If $(Z_n)_n \rightarrow Z$, pointwise, with all $Z_n$ in $\boldsymbol{U}$, then

(a) $(I(Z_n))_n \to I(Z)$, pointwise, with, for Y in **U,** $I(Y)$: $\theta \to I_\theta(Y)$.

(b) $(\mu(Z_n))_n \to \mu(Z)$, pointwise, with, for Y in **U,** $\mu(Y)$: $\theta \to \mu_\theta(Y)$.

(c) $(g(Z_n))_n \to g(Z)$, pointwise, with, for Y in **U,** $g(Y)$: $\theta \to g_\theta(Y)$.

(d) $(h(Z_n))_n \to h(Z)$, pointwise, with, for Y in **U,** $h(Y)$: $\theta \to h_\theta(Y)$.

(e) $(P(Z_n))_n \to P(Z)$, pointwise, with, for Y in **U,** $P(Y)$: $\theta \to P_\theta(Y) = Z(\theta)$.

(f) $(R(Z_n))_n \to R(Z)$, pointwise, with, for Y in **U,** $R(Y)$: $\theta \to R_\theta(Y)$.

(g) A similar result holds for (PED)-bundles.

We next show how parts (a), (b), and (c) can be strengthened. In this proof, we say that $\theta$ is admissible if $g_\theta(Z)$ and all $g_\theta(Z_n)$ exist.

Theorem 8.

If $(Z_n)_n \to Z$, pointwise, with all $Z_n$ in **U**, then

(a) $(I(Z_n))_n \to I(Z)$, uniformly in $\theta$.

(b) $(\mu(Z_n))_n \to \mu(Z)$, uniformly in $\theta$.

(c) if Z is continuous and $Z \neq$ **O**, then $(g(Z_n))_n \to g(Z)$, uniformly in $\theta$, with $\theta$ admissible.

(d) Point (c) is not valid if $Z =$ **O**.

Proof.

(a). As $|Z_n - Z|$ tends to zero and all functions $Z_n$ and $Z$ are decreasing, each function $Z_n$-$Z$ is integrable. The functions $Z_n$-$Z$ are moreover uniformly bounded, hence it follows from Theorem 2 (based on Arzelà's theorem) that $\int_0^T |Z_n(s) - Z(s)| ds \rightarrow 0$ .

As $\left| \int_0^\theta (Z_n(s) - Z(s)) ds \right| \leq \int_0^\theta |Z_n(s) - Z(s)| \, ds \leq \int_0^T |Z_n(s) - Z(s)| \, ds$ this proves that $\int_0^\theta Z_n(s) ds \rightarrow \int_0^\theta Z(s) ds$ uniformly in θ.

(b) For n fixed, $\mu(Z_n)$ is decreasing in the variable θ (because all $Z_n$ are decreasing). Moreover, it follows from point (a) that $\mu(Z_n)$ is pointwise decreasing to $\mu(Z)$ in the variable θ.  Now $\mu(Z)$ is differentiable, hence continuous in θ. Applying now the second Dini theorem on the sequence $f_n = \mu(Z_n)$ proves that $(\mu(Z_n))_n \rightarrow \mu(Z)$, uniformly in θ.

(c) Put $x = g_\theta(Z)$ then we know, see the proof of Theorem 3, that

$|g_\theta(Z_n) - g_\theta(Z)| \leq \frac{1}{\theta x} \left| \int_0^x Z_n(s) ds - \int_0^x Z(s) ds \right|$ . As $Z$ is decreasing, we know that $\theta x = \frac{1}{x} \int_0^x Z(s) ds \geq \mu_T(Z)$, hence: $|g_\theta(Z_n) - g_\theta(Z)| = \frac{1}{\mu_T(Z)} \left| \int_0^x Z_n(s) ds - \int_0^x Z(s) ds \right|$. Using part (1) of this theorem this proves part (2).

(d). Part (c) is not valid for $Z =$ **O**. We provide a counterexample. Let $Z_n$ be the constant function $a_n$ on [0,T] with $\lim_{n \to \infty} a_n = 0$. Then $(Z_n)_n$ tends uniformly, hence pointwise, on [0,T] to the function $Z =$ **O**. As, $g_\theta(Z_n) = h_\theta(Z_n) = a_n/\theta$, we see that this sequence tends to $g_\theta(Z) = h_\theta(Z) = 0$ pointwise, but not uniformly.

We next investigate what happens if $(Z_n)_n \to Z$, uniformly, with again all $Z_n$ , and hence $Z$, in **U**.

Theorem 9

If $(Z_n)_n \to Z$, uniformly, with all $Z_n$ in **U**, then

(a) $(I(Z_n))_n \to I(Z)$, uniformly in $\theta$.

(b) $(\mu(Z_n))_n \to \mu(Z)$, uniformly in $\theta$.

(c) If $Z \neq$ **O**, then also $(g(Z_n))_n \to g(Z)$, uniformly.

(d) If $\theta_0 = \inf\{\theta; \theta$ is admissible$\} > 0$, then $(h(Z_n))_n \to h(Z)$, uniformly. This property does not hold if $\theta_0 = 0$ (this happens e.g., when $Z(T) = 0$).

(e) $(P(Z_n))_n \to P(Z)$, uniformly.

(f) If $\theta_0 = \inf\{\theta; \theta$ is admissible$\} > 0$ then $(R(Z_n))_n \to R(Z)$, uniformly. This property does not hold if $\theta_0 = 0$ (this happens e.g., if $Z(T) = 0$).

Proof.

Points (a), (b), and (c) are just special cases of theorem 8 because uniform convergence implies pointwise convergence.

(d). From the proof of Theorem 4, with $f_m(s) = \theta s$, $x = h_\theta(Z)$ and $x_n = h_\theta(Z_n)$ we know:

$\theta |x - x_n| \leq |Z(x) - Z_n(x)|$.

Hence: $|h_\theta(Z) - h_\theta(Z_n)| \leq \frac{1}{\theta}|Z(x) - Z_n(x)|$

This shows that the uniform convergence of $(h(Z_n))_n$ follows from the uniform convergence of $(Z_n)_n$ and the fact that $\theta_0 > 0$.

(e). This is trivial as $|P_\theta(Z_n) - P_\theta(Z)| = |Z_n(\theta) - Z(\theta)|$

(f). From Theorem 6 we know that $\left|R_\theta^2(Z_n) - R_\theta^2(Z)\right|$

$$\leq \left|\int_0^{h_\theta(Z)} Z_n(s)ds - \int_0^{h_\theta(Z)} Z(s)ds\right| + \left|\int_{h_\theta(Z)}^{h_\theta(Z_n)} Z_n(s)ds\right|$$

Now, when n increases, the first term in this sum converges uniformly to zero because of the part (a) – using $h_\theta(Z)$ as $\theta$ - while the second term is smaller than $\max_n(Z_n(0), Z(0)) * |h_\theta(Z_n) - h_\theta(Z)|$. This term converges uniformly to zero because of point (d).

Remarks

1. Property (c) does not hold for Z = **O**. This is shown in part (d) of Theorem 8.

2. If $\theta_0 = \inf\{\theta; \theta$ is admissible$\} > 0$, then $(h(Z_n))_n \to h(Z)$, uniformly. This property does not hold if $\theta_0 = 0$ (this happens e.g. when Z(T)=0).

Let T and S be strict positive constants, and consider the sequence $(S_n)_{n=3, 4, ...}$ with $S_n = S/n$. Define $Z_n$ on [0,T] as the function whose graph linearly connects the points (0,S) and (T/2,S/2), is equal to S/n on the interval [3T/4, T], and which in between, linearly connects the points (T/2,S/2) and (3T/4, S/n). Then $(Z_n)_n \to Z$ uniformly with Z the function that coincides with all $Z_n$ on [0,T/2], then linearly connects the points (T/2,S/2) and (3T/4,0) and which is equal to zero on the interval [3T/4,T], see Fig.1

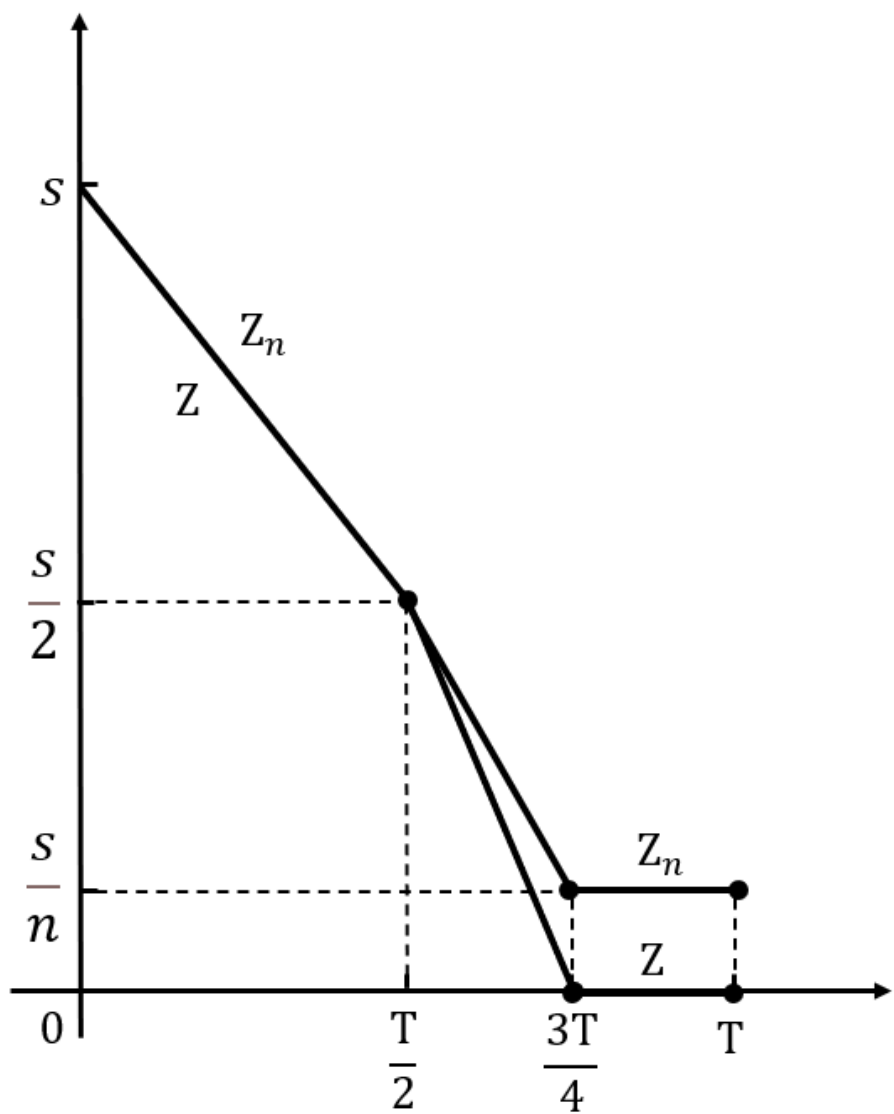

Fig. 1 Graphs illustrating Remark 2

We see that $Z \neq \mathbf{0}$, but $Z(T) = 0$ and $\theta_0 = 0$. From Theorem 7, we already know that $(h(Z_n))_n \to h(Z)$, pointwise. For $\theta < (4S)/3nT$ we have $h_\theta(Z_n) = S/n\theta$. Now, $h_\theta(Z_n)$ can be made smaller than a given ε > 0, by taking n > S/(εθ) = $n_0$. Hence the convergence is not uniform in θ.

3. If $θ_0$ = inf{θ; θ is admissible} > 0, then $(R(Z_n))_n \to R(Z)$, uniformly. This property does not hold if $θ_0$ = 0 (this happens e.g., when Z(T)=0).

We use the same example as in remark 2, with θ admissible and n ≤ m. Then

$$\left| R_\theta^2(Z_n) - R_\theta^2(Z_m) \right| = \left| \int_0^{h_\theta(Z_n)} Z_n(s)ds - \int_0^{h_\theta(Z_m)} Z_m(s)ds \right|$$

$$> \int_{3t/4}^{h_\theta(Z_n)} Z_n(s)ds - \int_{\frac{3T}{4}}^{h_\theta(Z_m)} Z_m(s)ds = \frac{S}{n}\left(\frac{S/n}{\theta} - \frac{3T}{4}\right) - \frac{S}{m}\left(\frac{S/m}{\theta} - \frac{3T}{4}\right)$$

This last expression tends to zero when n and m tend to ∞, but this does not happen uniformly in θ.

**A classification of bundles**

Based on our results related to the convergence of measures we will now come to a classification of bundles. We define three sets of bundles, denoted as (PC), (PC*), and (UC).

Definitions

A bundle, with bundle measure m belongs to the set (PC) iff

$[(Z_n)_n \to Z$, pointwise on [0,T]] => [$(m(Z_n))_n \to m(Z)$, pointwise in θ], where, of course, we only consider admissible values of θ.

A bundle, with bundle measure m belongs to the set (PC*) iff

$[(Z_n)_n \to Z$, pointwise on $[0,T]$, with $Z$ continuous$] => [(m(Z_n))_n \to m(Z)$, pointwise in $\theta]$, where, of course, we only consider admissible values of $\theta$. The difference between (PC) and (PC*) is that in (PC) we do not require that the limiting function $Z$ is continuous. Obviously, we have that (PC) $\subset$ (PC*).

A bundle, with bundle measure m belongs to (UC) iff

$[(Z_n)_n \to Z$, uniform on $[0,T]] => [(m(Z_n))_n \to m(Z)$, uniformly in $\theta]$, where, again, we only consider admissible values of $\theta$.

The next theorem follows directly from our earlier results. Note that a bundle is denoted by its bundle measure.

Theorem 10

(a) The bundles I, μ, g, h, $h^{(p)}$, R, and all (PED)-bundles belong to (PC).

(b) The bundles I, μ and P belong to (UC).

To come to a full description, we intend to show the following three results:

1) (UC) $\subset$ (PC*).

2) There exist bundles in (UC) that do not belong to (PC).

3) There exist bundles in (PC*) that do not belong to (PC) $\cup$ (UC).

Assuming that these results are shown we arrive at figure 2.

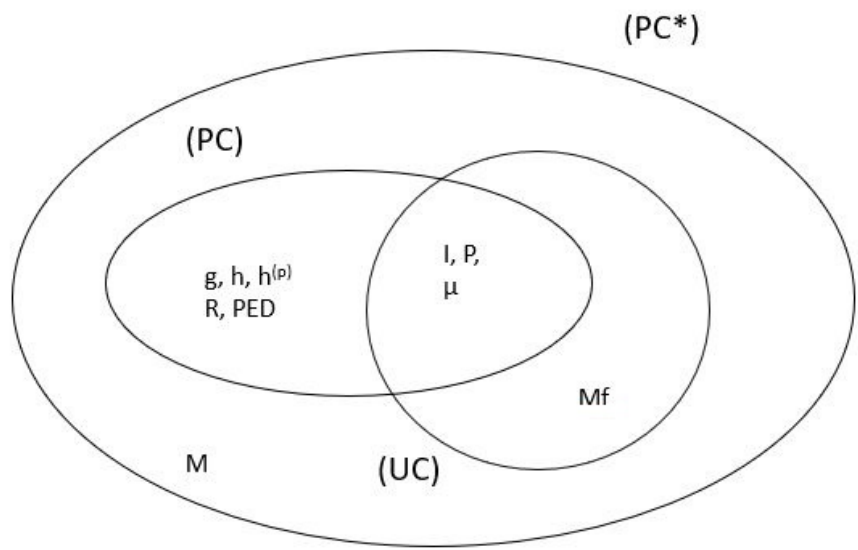

Fig. 2. A bundle classification (for Mf, see Theorem 12; for M
see Theorem 13)

Theorem 11. $(UC) \subset (PC^*)$

Proof. Assume that m belongs to (UC). If now $(Z_n)_n \to Z$,
pointwise on [0,T] with Z continuous, we have by Dini's second
theorem (recall that all $Z_n$ are decreasing) that $(Z_n)_n \to Z$,
uniformly. Because m belongs to (UC), we have then that
$(m(Z_n))_n \to m(Z)$, uniformly in θ (admissible). Hence, also
$(m(Z_n))_n \to m(Z)$, pointwise in θ, showing that m belongs to
$(PC^*)$.

Theorem 12. $(UC) \not\subset (PC)$

Proof. We provide an example of a bundle (measure) that
belongs to (UC) and does not belong to (PC). We already know
that if $\boldsymbol{Z} = \{Z_n, n \in \boldsymbol{N}\}$, with $Z_n(x) = 1 - x^n$, $x \in [0,1]$ , then

$(Z_n)_n$ is pointwise convergent to $Z(x)$, with $Z(x) = 1$ on $[0,1[$ and $Z(1) = 0$. Yet, this convergence is not uniform on $[0,1]$. This property also holds for every subsequence of $(Z_n)_n$. Hence there does not exist a uniform convergent sequence of functions in **Z**, and hence any bundle measure m on **Z** belongs to (UC).

Let now $\theta_0 \in ]0,1[$ and let f be a function that is not continuous in the point $1 = Z(\theta_0)$. Now, for Y in **Z** and $\beta \in [0,1]$, define

$$Mf_\beta(Y) = f(Y(\beta))$$

Then, $f(Z_n(\theta_0)) = f(1 - \theta_0^n)$ does not tend to $f(Z(\theta_0)) = f(1)$. Now, $Mf_{\theta_0}(Z_n) = f(Z_n(\theta_0))$ does not tend to $Mf_{\theta_0}(Z) = f(Z(\theta_0)) = f(1)$, showing that $Mf(Z_n)$ does not tend to $Mf(Z)$ pointwise.

Theorem 13. There exist bundles in (PC*) that do not belong to (PC) ∪ (UC).

Proof. We provide one example. For all functions $Y : [0,T] \to \mathbf{R}^+$ and all $\theta \in ]0,T]$ we define the bundle $M_\theta(Y) = \frac{1}{\theta}(\lim_{x \to \theta} Y(x))$. We first show that $M \in (PC*)$. Assume that $(Z_n)_n \to Z$, pointwise on $[0,T]$ with Z and all $Z_n$ continuous. Then we have for each $\theta \in [0,T]$: $M_\theta(Z_n) = \frac{1}{\theta}(\lim_{x \to \theta} Z_n(x)) = \frac{Z_n(\theta)}{\theta} \to \frac{Z(\theta)}{\theta} = \frac{\lim_{x \to \theta} Z(x)}{\theta} = M_\theta(Z)$. This shows that $M_\theta(Z_n) \to M_\theta(Z)$, pointwise. Hence $M \in (PC*)$.

Assume now that $(Z_n)_n \to Z$, pointwise on $[0,T]$ with all $Z_n$ continuous, but $Z: [0,T] \to \mathbf{R}^+$ not continuous in a point $\theta_0$ $\in\ ]0,T[$.     Then     $M_{\theta_0}(Z_n) = \lim_{x\to\theta_0} Z_n(x) = Z_n(\theta_0) \to Z(\theta_0) \neq \lim_{x\to\theta_0} Z(x) = M_{\theta_0}(Z)$ .

This shows that $M_\theta(Z_n)$ does not tend to $M_\theta(Z)$ pointwise for all $\theta \in ]0,T[$. Hence $M \notin$ (PC).

We still have to show that $M \notin$ (UC).

Let $(a_n)_n$ be a sequence of positive real numbers converging to the positive real number a. Let $(Z_n)_n$ be a sequence of constant functions on [0,T] with $Z_n(x) = a_n$. Then $(Z_n)_n \to Z$, uniformly with $Z(x) = a$. Now $M_\theta(Z_n) = a_n/\theta$ which converges to $a/\theta$, but this convergence is not uniformly.

This result concludes the explanation of Fig.2

## Conclusion

This paper studied the issue of stability of impact measures and bundles, through convergence properties. It is proposed to include these properties in the study of impact. This is an aspect that - to the best of our knowledge - has not been addressed so far in the informetric literature. We showed that pointwise convergence is maintained by all well-known impact bundles (such as the h-, g-, and R-bundle) and that the μ-bundle (and variants) even maintain uniform convergence. Based on these results, a classification of impact bundles is given.

Acknowledgment. The author thanks Li Li (National Science Library, CAS) for drawing Figure 1.